\magnification1200
\input amssym.def
\input amssym.tex

\hsize=13.5truecm
\baselineskip=16truept plus.3pt minus .3pt

\font\secbf=cmb10 scaled 1200
\font\eightrm=cmr8
\font\sixrm=cmr6

\font\eighti=cmmi8

\font\sixi=cmmi6
\skewchar\eighti='177 \skewchar\sixi='177

\font\eightsy=cmsy8
\font\sixsy=cmsy6
\skewchar\eightsy='60 \skewchar\sixsy='60

\font\eightit=cmti8

\font\eightbf=cmbx8
\font\sixbf=cmbx6

\let\sc=\tensc

\font\eightsc=cmcsc10 scaled 800
%\font\cyr=wncyr10
%\font\cyrbf=wncyb10
%\font\cyrit=wncyi10
%\font\abscyr=wncyr8             % roman do streszcze¤ i literatury
%\font\itcyr=wncyi8              % italic do literatury
%\font\sccyr=wncysc10 scaled 800
    % Font na tytu'y
 % Font na podtytu'y
\font\secbf=cmb10 scaled 1200
\font\subsecfont=cmb10 scaled \magstephalf
%%%%%%%%%%%%%%%%%%%%%%%%%%%%%%%%%%%%%%%%%%%%%%%%%%
\font\amb=cmmib10

\font\ambi=cmmib10 scaled 700

\newfam\mbfam 

\textfont\mbfam\amb \scriptfont\mbfam\ambi

%%%%%%%%%%%%%%%%%%%%%%%%%%%%%%%%%%%%%%%%%%%%%%%%%%%

\def\aa{\def\rm{\fam0\eightrm}%
  \textfont0=\eightrm \scriptfont0=\sixrm \scriptscriptfont0=\fiverm
  \textfont1=\eighti \scriptfont1=\sixi \scriptscriptfont1=\fivei
  \textfont2=\eightsy \scriptfont2=\sixsy \scriptscriptfont2=\fivesy
  \textfont3=\tenex \scriptfont3=\tenex \scriptscriptfont3=\tenex
  \def\sc{\eightsc}
  \def\it{\fam\itfam\eightit}%
  \textfont\itfam=\eightit
  \def\bf{\fam\bffam\eightbf}%
  \textfont\bffam=\eightbf \scriptfont\bffam=\sixbf
   \scriptscriptfont\bffam=\fivebf
  \normalbaselineskip=9.7pt
  \setbox\strutbox=\hbox{\vrule height7pt depth2.6pt width0pt}%
  \normalbaselines\rm}

\def\Proof{\vskip12pt\noindent{\bf Proof.} }

\def\Remark#1{\vskip12pt\noindent{\bf Remark #1}}

\def\m@th{\mathsurround=0pt}

\def\cc#1{\hbox to .89\hsize{$\displaystyle\hfil{#1}\hfil$}\cr}
\def\lc#1{\hbox to .89\hsize{$\displaystyle{#1}\hfill$}\cr}
\def\rc#1{\hbox to .89\hsize{$\displaystyle\hfill{#1}$}\cr}

\def\eqal#1{\null\,\vcenter{\openup\jot\m@th
  \ialign{\strut\hfil$\displaystyle{##}$&&$\displaystyle{{}##}$\hfil
      \crcr#1\crcr}}\,}

\def\section#1{\vskip 22pt plus6pt minus2pt\penalty-400
        {{\secbf
        \noindent#1\rightskip=0pt plus 1fill\par}}
        \par\vskip 12pt plus5pt minus 2pt
        \penalty 1000}

\def\subsection#1{\vskip 20pt plus6pt minus2pt\penalty-400
        {{\subsecfont
        \noindent#1\rightskip=0pt plus 1fill\par}}
        \par\vskip 8pt plus5pt minus 2pt
        \penalty 1000}

\def\subsubsection#1{\vskip 18pt plus6pt minus2pt\penalty-400
        {{\subsecfont
        \noindent#1}}
        \par\vskip 7pt plus5pt minus 2pt
        \penalty 1000}

\def\center#1{{\begingroup \leftskip=0pt plus 1fil\rightskip=\leftskip
\parfillskip=0pt \spaceskip=.3333em \xspaceskip=.5em \pretolerance 9999
\tolerance 9999 \parindent 0pt \hyphenpenalty 9999 \exhyphenpenalty 9999
\par #1\par\endgroup}}

\def\\{\hfill\break}

\def\kwadrat{\null\ \hfill\null\ \hfill$\square$}
\def\mida#1{{{\null\kern-4.2pt\left\bracevert\vbox to 6pt{}\!\hbox{$#1$}\!\right\bracevert\!\!}}}
\def\midy#1{{{\null\kern-4.2pt\left\bracevert\!\!\hbox{$\scriptstyle{#1}$}\!\!\right\bracevert\!\!}}}

\def\diagin{\hbox{---}\hskip-11.5pt\intop}
\def\diagint{\hbox{--}\hskip-8.2pt\intop}
\def\diagintop{\mathop{\mathchoice
{{\diagin}}%
{{\diagint}}%
{{\diagint}}%
{{\diagint}}%
}\limits}

\def\divv{{{\rm div}\,}}
\def\rot{{\rm rot}\,}
\def\const{{\rm const}}

\def\today{${\scriptscriptstyle\number\day-\number\month-\number\year}$}
\footline={{\hfil\rm\the\pageno\hfil${\scriptscriptstyle\rm\jobname}$\ \ \today}}

%%%%%%%%%%%%%%%%%%%%%%%%%%%%%%%%%%
\def\ifnextchar#1#2#3{\bgroup
  \def\reserveda{\ifx\reservedc #1 \aftergroup\firstoftwo
    \else \aftergroup\secondoftwo\fi\egroup{#2}{#3}}%
  \futurelet\reservedc\ifnch
  } 
\def\ifnch{\ifx \reservedc \sptoken \expandafter\xifnch
      \else \expandafter\reserveda
      \fi} 
\def\firstoftwo#1#2{#1}
\def\secondoftwo#1#2{#2} 
\def\tempswafalse{\let\iftempswa\iffalse}
\def\tempswatrue{\let\iftempswa\iftrue} 

\def\cite{\ifnextchar [{\tempswatrue\citea}{\tempswafalse\citeb}}
\def\citea[#1]#2{[#2, #1]}
\def\citeb#1{[#1]}\phantom{]}
%%%%%%%%%%%%%%%%%%%%%%%%%%%%%%%%%%%%%%%%%

\def\\{\hfil\break}

\def\D{{\Bbb D}}
\def\N{{\Bbb N}}
\def\R{{\Bbb R}}
\def\T{{\Bbb T}}

\def\divv{{\rm div}\,}
\def\rot{{\rm rot}}
\def\meas{{\rm meas}}
\def\esssup{{\rm esssup}}

\center{\secbf Stability of two-dimensional solutions to the Navier-Stokes 
equations in cylindrical domains under Navier boundary\break conditions}

\vskip1cm
\centerline{\bf W. M. Zaj\c aczkowski}
\vskip1cm
%\item{$^1$}
\noindent Institute of Mathematics Polish Academy of Sciences,\\ 
\'Sniadeckich 8, 00-656 Warsaw, Poland,\\
e-mail:wz@impan.pl;\\
Institute of Mathematics and Cryptology, Cybernetics Faculty,\\ 
Military University of Technology,\\
Kaliskiego 2, 00-908 Warsaw, Poland
\vskip1cm

\noindent
{\bf Abstract.} 
The Navier-Stokes motions in a cylindrical domain with Navier boundary 
conditions are considered. First the existence of global regular 
two-dimensional solutions is proved. The solutions are such that bounded with 
respect to time norms are controlled by the same constant for all $t\in\R_+$. 
Assuming that the initial velocity and the external force are sufficiently 
close to the initial velocity and the external force of the two-dimensional 
solutions we prove existence of global three-dimensional solutions which 
remain close to the two-dimensional solutions for all time. In this way we 
mean stability of two-dimensional solutions. Thanks to the Navier boundary 
conditions the nonlinear term in two-dimensional Navier-Stokes equations does 
not have any influence on the form of the energy estimate. This implies that 
stability is proved without any structural restrictions on the external force, 
initial data and viscosity.


\vfil\eject 

\section{1. Introduction}

In this paper we prove stability of two-dimensional solutions in a set of 
three-dimensional motions of the Navier-Stokes equations in cylindrical domain 
$D=\Omega\times(-a,a)$, where $\Omega\subset\R^2$ and $L=2a$ is the length of 
the cylinder. The tree-dimensional motions satisfy the following 
initial-boundary value problem
$$\eqal{
&v_t+v\cdot\nabla v+\nu\rot\rot v+\nabla p=f\quad &{\rm in}\ \ 
D_+\equiv D\times\R_+,\cr
&\divv v=0\quad &{\rm in}\ \ D_+,\cr
&v\cdot\bar n=0\quad &{\rm on}\ \ S_+=S\times\R_+,\cr
&\bar n\times\rot v=0\quad &{\rm on}\ \ S_+,\cr
&v|_{t=0}=v(0)\quad &{\rm in}\ \ D,\cr}
\leqno(1.1)
$$
where $v=(v_1(x,t),v_2(x,t),v_3(x,t))\in\R^3$ is the velocity of the fluid, 
$p=p(x,t)\in\R$ is the pressure, $f=(f_1(x,t),f_2(x,t),f_3(x,t))\in\R^3$ is 
the external force field. By $x=(x_1,x_2,x_3)$ are denoted the Cartesian 
coordinates such that $x_3$-axis is parallel to the cylinder and is located 
inside it. By the dot we denote the scalar product in $\R^3$.

\noindent
The Navier-Stokes equations $(1.1)_1$ follow from the formula
$$
-\Delta v=\rot\rot v
\leqno(1.2)
$$
which holds for divergence free vectors $v$. Finally, by $\nu$ we denote the 
positive viscosity coefficient.

Moreover, $\bar n$ is the unit outward vector normal to $S$. The boundary $S$ 
is split into two parts, $S=S_1\cup S_2$, where $S_1$ is parallel to the 
$x_3$-axis and $S_2$ is perpendicular. Additionally $S_2=S_2(-a)\cup S_2(a)$, 
where $S_2(b)$ meets $x_3$-axis at $x_3=b$.

Our aim is to prove existence of global regular nonvanishing with time 
solutions to problem (1.1). For this we need that the external force does not 
converge to zero as time goes to infinity. Then to reduce restrictions on the 
external force we introduce the quantities
$$
p'=x\cdot\diagintop fdx,\quad \bar p=p-p',\quad \bar f=f-\diagintop fdx,
\leqno(1.3)
$$
where
$$
\diagintop fdx={1\over|D|}\intop_Df(x)dx\quad {\rm and}\quad |D|=\meas D.
$$
Then problem (1.1) takes the form
$$\eqal{
&v_t+v\cdot\nabla v+\nu\rot\rot v+\nabla\bar p=\bar f,\cr
&\divv v=0,\cr
&v\cdot\bar n=0,\quad \bar n\times\rot v=0,\cr
&v|_{t=0}=v(0).\cr}
\leqno(1.4)
$$
Up to now there is not possible to prove existence of global regular solutions 
to problem (1.1). Therefore we restrict our considerations to show existence 
of global regular solutions which remain sufficiently close to two-dimensional 
global correspondingly regular solutions. It is well known that such 
two-dimensional solutions exist. Since we need a special behavior of 
two-dimensional solutions we show their existence in Section 3.

By two-dimensional motions we mean such solutions to (1.1) that 
$v=w=(w_1(x_1,x_2,t),w_2(x_1,x_2,t))\in\R^2$, $p=\eta(x_1,x_2,t)\in\R$ and 
$f=h=(h_1(x_1,x_2,t),h_2(x_1,x_2,t))\in\R^2$. Hence the two-dimensional 
motions satisfy
$$\eqal{
&w_t+w\cdot\nabla w+\nu\tilde\rot\rot^{(2)}w+\nabla\eta=h\quad &{\rm in}\ \ 
\Omega\times\R_+\equiv\Omega_+,\cr
&\divv w=0\quad &{\rm in}\ \ \Omega_+\cr
&w\cdot\bar n=0\quad &{\rm on}\ \ S_0\times\R\equiv S_{0+},\cr
&\rot^{(2)}w=0\quad &{\rm on}\ \ S_{0+},\cr
&w|_{t=0}=w(0)\quad &{\rm in}\ \ \Omega,\cr}
\leqno(1.5)
$$
where $S_0=\partial\Omega$ and $\rot^{(2)}w=w_{2,x_1}-w_{1,x_2}$, 
$\tilde\rot\varphi=(\varphi_{,x_2}-\varphi_{,x_1})$. 
Comparing to (1.1) we see that $S_1=S_0\times(-a,a)$.

To examine problem (1.5) we need transformation of type (1.3) applied to the 
two-dimensional case. Therefore, in this case we introduce
$$
\eta'=x'\cdot\diagintop hdx',\quad \bar\eta=\eta-\eta',\quad 
\bar h=h-\diagintop hdx',
\leqno(1.6)
$$
where $x'=(x_1,x_2)$, $\diagintop hdx'={1\over|\Omega|}\intop_\Omega h(x')dx'$.

\noindent
Using (1.6) problem (1.5) takes the form
$$\eqal{
&w_t+w\cdot\nabla w+\nu\tilde\rot\rot^{(2)}w+\nabla\bar\eta=\bar h\quad 
&{\rm in}\ \ \Omega_+,\cr
&\divv w=0\quad &{\rm in}\ \ \Omega_+,\cr
&w\cdot\bar n=0,\ \ \rot^{(2)}w=0\quad &{\rm on}\ \ S_{0+},\cr
&w|_{t=0}=w(0)\quad &{\rm in}\ \ \Omega.\cr}
\leqno(1.7)
$$
Since we consider incompressible motions we can assume without any 
restrictions that $f$ and $h$ are divergence free.

\noindent
To show stability of two-dimensional solutions we introduce the quantities
$$
u=v-w,\quad q=\bar p-\bar\eta,\quad g=\bar f-\bar h
\leqno(1.8)
$$
which are solutions to the problem
$$\eqal{
&u_t+u\cdot\nabla u=\nu\rot^2u+\nabla q=-w\cdot\nabla u-u\cdot\nabla w+g\quad 
&{\rm in}\ \ D_+,\cr
&\divv u=0\quad &{\rm in}\ \ D_+,\cr
&u\cdot\bar n=0,\ \ \bar n\times\rot u=0\quad &{\rm on}\ \ S_+,\cr
&u|_{t=0}=u(0)\quad &{\rm in}\ \ D.\cr}
\leqno(1.9)
$$

\Remark{1.1.} 
The operator $\rot^2$ and the boundary conditions $(1.9)_3$ hold for function 
$v$. Now we show that they are also satisfied for a solution $w$ to problem 
(1.7). Using that $w=(w_1,w_2,0)$ and $\divv w=0$ we have
$$\eqal{
&\rot^2w=-\Delta^{(2)}w=-\left(\matrix{
w_{1,x_1x_1}+w_{1,x_2x_2}\cr w_{2,x_1x_1}+w_{2,x_2x_2}\cr}\right)\cr
&=\left(\matrix{(w_{1,x_1}+w_{2,x_2})_{,x_1}-(w_{1,x_1x_1}+w_{1,x_2x_2})\cr
(w_{1,x_1}+w_{2,x_2})_{,x_2}-(w_{2,x_1x_1}+w_{2,x_1x_2})\cr}\right)\!=\!
\left(\matrix{\phantom{-}(w_{2,x_1}-w_{1,x_2})_{,x_2}\cr 
-(w_{2,x_1}-w_{1,x_2})_{,x_1}\cr}\right)\cr
&=\tilde\rot\rot^{(2)}w,\cr}
$$
where $\Delta^{(2)}=\partial_{x_1}^2+\partial_{x_2}^2$ and $\rot^2=\rot\rot$. 
Therefore, formulation of the operator $\rot^2$ in $(1.9)_1$ is right.

\noindent
To satisfy boundary conditions $(1.9)_3$ we have to introduce the tangent and 
normal vectors to $S_1$ and $S_2$. From the geometry of a cylinder we have 
$\bar n|_{S_1}=\bar n|_{S_0}$. Let $S_0$ be described by a sufficiently 
regular function $\varphi(x_1,x_2)=0$. Then 
$\bar n|_{S_0}={\nabla\varphi\over|\nabla\varphi|}$. Next the tangent vector 
to $S_0$, denoted by $\bar\tau$, equals 
$\bar\tau={\bar\nabla\varphi\over|\nabla\varphi|}$, where 
$\bar\nabla\varphi=(-\varphi_{,x_2},\varphi_{,x_1})$. We have two tangent 
vectors to $S_1$: $\bar\tau_1=(\bar\tau,0)$, $\bar\tau_2=(0,0,1)$. On $S_2$ 
we have $\bar n|_{S_2}=(0,0,1)$, $\bar\tau_1|_{S_2}=(1,0,0)$, 
$\bar\tau_2|_{S_2}=(0,1,0)$.

\noindent
Since $\bar n|_{S_1}=\bar n|_{S_0}$ we have that $u\cdot\bar n=0$ on $S_1$. 
We have also that $w\cdot\bar n|_{S_2}=0$, so finally $u\cdot\bar n|_S=0$ 
holds. Next we examine the condition
$$
\bar n\times\rot w|_S=0
\leqno(1.10)
$$
On $S_1$ it is equivalent to
$$\eqal{
&\bar\tau_1\cdot\rot w|_{S_1}=0,\cr
&\bar\tau_2\cdot\rot w|_{S_1}=0,\cr}
\leqno(1.11)
$$
where the second condition $(1.11)_2$ equals $w_{1,x_2}-w_{2,x_1}|_{S_1}=0$ 
so\break $\rot^{(2)}w|_{S_1}=0$ in view of boundary condition $(1.5)_4$. To satisfy 
$(1.11)_1$ we express it explicitly in the form
$$
\tau_1(w_{2,x_3}-w_{3,x_2})+\tau_2(w_{3,x_1}-w_{1,x_3})|_{S_1}=0.
$$
It holds because $w_{i,x_3}=0$, $i=1,2$, and $w_3=0$.

\noindent
Similarly, we show that $\bar n\times\rot w|_{S_2}=0$. Hence (1.10) holds.

\noindent
Now we present results of this paper. The introduced norms in these 
formulations are defined in Section 2. Remark 3.3 yields

\proclaim Theorem 1. 
Let $T>0$ be given. Let $w(0)\in H^1(\Omega)$, 
$\bar h\in L_2(kT,(k+1)T;L_2(\Omega))$ for any $k\in\N_0$. Then there exists 
a solution to problem (1.7) such that $w\in V_2^1(kT,(k+1)T;\Omega)$ and
$$
\|w(kT)\|_{H^1(\Omega)}\le c\bar A_1,
$$
$$
\|w(t)\|_{H^1(\Omega)}^2+\intop_{kT}^t\|w(t')\|_{H^2(\Omega)}^2dt'\le 
c\bar A_1^2,
$$
where $t\in(kT,(k+1)T)$, $k\in\N_0$ and
$$
\bar A_1^2=\sup_k{1\over\nu}\intop_{kT}^{(k+1)T}\|\bar h(t)\|_{L_2(\Omega)}^2dt
+\|w(0)\|_{H^1(\Omega)}^2.
$$

\noindent
Lemma 3.4 implies

\proclaim Theorem 2. \hskip-1pt
Let \hskip-1pt the \hskip-1pt assumptions \hskip-1pt of \hskip-1pt Theorem \hskip-1pt 1 \hskip-1pt hold. \hskip-1pt Let \hskip-1pt $w(0)\in B_{\sigma,2}^1(\Omega)$,\break 
$\bar h\in L_2(kT,(k+1)T;L_\sigma(\Omega))$, $\sigma>3$, $k\in\N_0$. Then
$$
\|w(t)\|_{W_\sigma^1(\Omega)}+
\|w\|_{W_{\sigma,2}^{2,1}(\Omega\times(kT,(k+1)T))}\le\bar A_2,
$$
where $\bar A_2$ depends on $\bar A_1$, 
$$
\|w(0)\|_{B_{\sigma,2}^1(\Omega)},\bigg(\sup_{k\in\N_0}\intop_{kT}^{(k+1)T}
\|\bar h(t)\|_{L_\sigma(\Omega)}^2dt\bigg)^{1/2}.
$$

\noindent
From Lemma 4.2 we have

\proclaim Theorem 3. 
Let the assumptions of Theorems 1 and 2 hold. Let $c_*\in(0,\nu]$. Let 
$\gamma_*$ be so small that ${c_*\over2}\le\nu-{c_0\over\nu^3}\gamma_*^4$, 
where $c_0$ is the constant from (4.20). Let $\gamma\le\gamma_*$. Let
$$\eqal{
&\|u(0)\|_{H^1(\Omega)}^2\le\gamma,\cr
&{c\over\nu}\|g(t)\|_{L_2(D)}^2\le{c_*\over4}\gamma.\cr}
$$
Then
$$
\|u(t)\|_{H^1(\Omega)}^2\le\gamma\quad {\sl for}\quad t\in\R_+.
$$

From Theorems 1 and 3 and Remarks 3.5 and 4.3 we have

\proclaim Theorem 4. 
Assume that $f\in L_2(kT,(k+1)T;L_2(\Omega))$, $v(0)\in H^1(\Omega)$, 
$k\in\N_0$. Then there exists a global regular solution to (1.1) such that 
$v=w+u$, $p=\eta+q$, $f=h+g$ and $v\in W_2^{2,1}(\Omega\times(kT,(k+1)T))$, $\nabla p\in L_2(\Omega\times(kT,(k+1)T))$.

In \cite{22} problem (1.1) is considered in the periodic box. In this case global existence of two-dimensional solutions such that 
$w\in V_2^1(kT,(k+1)T;\Omega)$ is proved under very restrictive relation between $T$, $\nu$, $\bar A_1$. The relation holds for sufficiently large $T$, $\nu$ and correspondingly small $\bar A_1$. In this paper we omit the restriction by considering problem (1.1) in a cylindrical domain with the Navier boundary condition. This, in view of (2.16), gives that the estimate
$$
\|w\|_{W_2^1(kT,(k+1)T;\Omega)}\le c_0,\quad k\in\N_0,
\leqno(1.12)
$$
holds without any restrictions. The main aim in \cite{22} and in this paper is to show that constant $c_0$ in (1.12) does not depend on $k$.

\noindent
This guarantees that the two-dimensional solutions does not increase in time. Hence, also stability of two-dimensional solutions can be proved.

The first results connected with the stability of global regular solutions to 
the nonstationary Navier-Stokes equations were proved by Beirao da Veiga and 
Secchi \cite{3}, followed by Ponce, Racke, Sideris and Titi \cite{18}. 
Paper \cite{3} is concerned with the stability in $L_p$-norm of a strong 
three-dimensional solution of the Navier-Stokes system with zero external 
force in the whole space. In \cite{18}, assuming that the external force is 
zero and a three-dimensional initial function is close to 
a two-dimensional one in $H^1(\R^3)$, the authors showed the existence of 
a global strong solution in $\R^3$ which remains close to a two-dimensional 
strong solution for all times. In \cite{17} Mucha obtained a similar result 
under weaker assumptions about the smallness of the initial velocity 
perturbation.

In the class of weak Leray-Hopf solutions the first stability result was 
obtained by Gallagher \cite{8}. She proved the stability of 
two-dimensional solutions of the Navier-Stokes equations with 
periodic boundary conditions under three-dimensional perturbations both in 
$L_2$ and $H^{1\over2}$ norms.

The stability of nontrivial periodic regular solutions to the Navier-Stokes 
equations was studied by Iftimie \cite{10} and by Mucha \cite{15}. 
The paper \cite{15} is devoted to the case when the external force is 
a potential belonging to $L_{r,loc}(\T^3\times[0,\infty))$ and when the intial 
data belongs to the space $W_r^{2-2/r}(\T^3)\cap L_2(\T^3)$, where $r\ge2$ and 
$\T$ is a torus. Under the assumption that there exists a global solution with 
data of regularity mentioned above and assuming that small perturbations 
of data have the same regularity as above, the author proves that 
perturbations of the velocity and the gradient of the pressure remain small in 
the spaces $W_r^{2,1}(\T^3\times(k,k+1))$ and $L_r(\T^3\times(k,k+1))$, 
$k\in\N$, respectively. Paper \cite{10} contains results concerning the 
stability of two-dimensional regular solutions to the Navier-Stokes system in 
a three-dimensional torus but here the initial data in the three-dimensional 
problem belong to an anisotropic space of functions having different 
regularity in the first two directions than in the third direction, and the 
external force vanishes. Moreover, Mucha \cite{16} studies the stability of 
regular solutions to the nonstationary Navier-Stokes system in $\R^3$ assuming 
that they tend in $W_r^{2,1}$ spaces $(r\ge2)$ to constant flows.

The papers of Auscher, Dubois and Tchamitchian \cite{1} and of Gallagher, 
Iftimie and Planchon \cite{9} concern the stability of global regular 
solutions to the Navier-Stokes equations in the whole space $\R^3$ with zero 
external force. These authors show that the norms of the considered 
solutions decay as $t\to\infty$.

It is worth mentioning the paper of Zhou \cite{25}, who proved the asymptotic 
stability of weak solutions $u$ with the property: 
$u\in L_2(0,\infty,BMO)$ to the Navier-Stokes equations in $\R^n$, $n\ge3$, 
with a force vanishing as $t\to\infty$.

An interesting result was obtained by Karch and Pilarczyk \cite{11}, who 
concentrate on the stability of Landau solutions to the Navier-Stokes system 
in $\R^3$. Assuming that the external force is a singular distribution they 
prove the asymptotic stability of the solution under any $L_2$-perturbation.

Paper \cite{7} of Chemin and Gallagher is devoted to the stability 
of some unique global solution with large data in a very weak sense.

Finally, the stability of Leray-Hopf weak solutions has recently been examined 
by Bardos et al. \cite{2}, where equations with vanishing external force are 
considered. That paper concerns the following three cases: two-dimensional 
flows in infinite cylinders under three-dimensional perturbations which are 
periodic in the vertical direction; helical flows in circular cylinders under 
general three-dimensional perturbations; and axisymmetric flows under general 
three-dimensional perturbations. The theorem concerning the first case 
extends a result obtained by Gallagher \cite{8} for purely periodic 
boundary conditions.

\noindent
Most of the papers discussed above concern to the case with zero external 
force (\cite{1--3}, \cite{7--10}, \cite{17}, \cite{18}) 
or with force which decays as $t\to\infty$ (\cite{18}). Exceptions are 
\cite{11, 15, 16}, where very special external forces, which are singular 
distributions in \cite{11} or potentials in \cite{15--16}, are 
considered. However, the case of potential forces is easily reduced to the 
case of zero external forces.

The aim of our paper is to prove the stability result for a large class of 
external forces $f_s$ which do not produce solutions decaying as $t\to\infty$. 

It is essential that our stability results are obtained together with the 
existence of a global strong three-dimensional solution close to 
a two-dimensional one.

The paper is divided into two main parts. In the first we prove existence of 
global strong two-dimensional solutions not vanishing as $t\to\infty$ because 
the external force does not vanish either. To prove existence of such 
solutions we use the step by step method. For this purpose we have to show 
that the data in the time interval $[kT,(k+1)T]$, $k\in\N$, do not increase 
with $k$. We do not need any restrictions on the time step $T$.
In the second part we prove existence of three-dimensional solutions that 
remain close to two-dimensional solutions. For this we need the initial 
velocity and the external force to be sufficiently close in apropriate norms 
to the initial velocity and the external force of the two-dimensional 
problems.

The proofs of this paper are based on the energy method. Thanks to the Navier boundary conditions the nonlinear term in the two-dimensional Navier-Stokes equations does not have any influence on the 
form of the energy estimate. The proofs of global existence 
which follow from the step by step technique are possible thanks to the natural 
decay property of the Navier-Stokes equations. This is mainly used in the 
first part of the paper (Section 3). To prove stability (Section 4) we use 
smallness of data $(v(0)-v_s(0)),(f-f_s)$ and a contradiction argument applied 
to the nonlinear ordinary differential inequality (4.24). The paper is 
a generalization of results from \cite{22, 24}, where the periodic case is 
considered.

We restrict ourselves to prove estimates only,because existence follows easily 
by the Faedo-Galerkin method.

The paper is organized as follows. In Section 2 we introduce notation and give 
some auxiliary results. Section 3 is devoted to the existence of 
a two-dimensional solution. It also contains some useful estimates of the 
solution. In Section 4 we prove the existence of a global strong solution 
to problem (1.1) close to the two-dimensional solution for all time.

\section{2. Notation and auxiliary results}

Let $\N_0=\N\cup\{0\}$. By $L_p(\Omega)$, $p\in[1,\infty]$, 
$\Omega\subset\R^n$ we denote the Lebesgue space of integrable functions and 
by $H^s(\Omega)$, $s\in\N_0$, $\Omega\subset\R^n$, the Sobolev space of 
function with the finite norm
$$
\|u\|_{H^s}=\|u\|_{H^s(\Omega)}=\bigg(\sum_{|\alpha|\le s}\intop_\Omega
|D_x^\alpha u(x)|^2dx\bigg)^{1/2},
$$
where $D_x^\alpha=\partial_{x_1}^{\alpha_1}\dots\partial_{x_n}^{\alpha_n}$, 
$|\alpha|=\alpha_1+\alpha_2+\dots+\alpha_n$, $\alpha_i\in\N_0$, 
$i=1,\dots,n$, $n=2,3$.
Let $u=(u_1,\dots,u_n)$ be a vector. Then $|u|=\sqrt{u_1^2+\cdots+u_n^2}$.

Next we introduce the anisotropic Lebesgue and Sobolev spaces with the mixed 
norms. $L_{p_1,p_2}(\Omega\times(0,T))$ and 
$W_{p_1,p_2}^{s,s/2}(\Omega\times(0,T))$, $p_1,p_2\in(1,\infty)$, 
${s\over2}\in\N$, are spaces with the following finite norms
$$\eqal{
&\|u\|_{L_{p_2}(0,T;L_{p_1}(\Omega))}\equiv\|u\|_{L_{p_1,p_2}(\Omega^T)}\!=\!\!
\bigg(\intop_0^T\!\bigg(\!\intop_\Omega|u(x,t)|^{p_1}dx\!\bigg)^{p_2/p_1}\!\!\!\!
dt\bigg)^{1/p_1}\!\!\!\!\!\!<\infty,\cr
&\|u\|_{W_{p_1,p_2}^{s,s/2}(\Omega^T)}=\!\!\sum_{|\alpha|+2\alpha_0\le s}
\bigg(\intop_0^T\bigg(\intop_\Omega
|D_x^\alpha\partial_t^{\alpha_0}u(x,t)|^{p_1}dx\bigg)^{p_2/p_1}\!\!\!dt\bigg)^{1/p_1}\!\!\!
<\infty,\cr}
$$
where $s$ is even.

\noindent
From [4, Ch. 4, Sect. 18] we recall the definition of the Besov spaces used in 
this paper. By $B_{p,q}^l(\Omega)$, $l\in\R_+$, $p,q\in[1,\infty]$, $\Omega\subset\R^n$, we denote the linear normed space of functions $u=u(x)$, $x\in\Omega$, with the finite norm
$$
\|u\|_{B_{p,q}^l(\Omega)}=\|u\|_{L_p(\Omega)}+\sum_{i=1}^n\bigg(\intop_0^{h_0}
\bigg|{\|\Delta_i^m(h;\Omega)\partial_i^ku\|_{L_p(\Omega)}\over 
h^{l-k}}\bigg|^q{dh\over h}\bigg)^{1/q},
$$
where $m>l-k>0$, $m,k\in\N$, $\Delta_i(h)u(x)=u(x+h\bar e_i)-u(x)$, 
$\Delta_i^m(h)u(x)=\Delta_i(\Delta_i^{m-1}(h)u(x))$, $\bar e_i$ -- versor of 
the $i$-th axis, $i=1,\dots,n$.

\noindent
Moreover,
$$
\Delta_i^m(h;\Omega)u(x)=\cases{
\Delta_i^m(h;\Omega)u(x)&for $[x,x+mhe_i]\in\Omega$,\cr
0&for $[x,x+mhe_i]\not\in\Omega$.\cr}
$$

In this paper the energy method is used to show main estimates and existence. 
For this purpose we need space $V_2^k(T_1,T_2;\Omega)$, $k\in\N$, with the 
finite norm
$$
\|u\|_{V_2^k(T_1,T_2;\Omega)}=\bigg(\esssup_{T_1\le t\le T_2}
\|u(t)\|_{H^k(\Omega)}^2+\intop_{T_1}^{T_2}
\|u(t)\|_{H^{k+1}(\Omega)}^2dt\bigg)^{1/2}.
$$
Let us consider the problem
$$\eqal{
&w_{2,x_1}-w_{1,x_2}=b\quad &{\rm in}\ \ \Omega,\cr
&w_{1,x_1}+w_{2,x_2}=0\quad &{\rm in}\ \ \Omega,\cr
&w\cdot\bar n=0\quad &{\rm on}\ \ S_0,\cr}
\leqno(2.1)
$$

\proclaim Lemma 2.1. 
Assume that $b\in H^s(\Omega)$, $s\in\N_0$. Then there exists a solution to 
problem (2.1) such that $w\in H^{s+1}(\Omega)$ and
$$
\|w\|_{H^{s+1}(\Omega)}\le c\|b\|_{H^s(\Omega)}.
\leqno(2.2)
$$

\Proof 
Equation $(2.1)_2$ implies existence of potential $\psi$ such that 
$w_1=-\psi_{,x_2}$, $w_2=\psi_{,x_1}$. Then problem (2.1) takes the form
$$\eqal{
&\Delta^{(2)}\psi=b\quad &{\rm in}\ \ \Omega,\cr
&\bar\tau\cdot\nabla\psi=0\quad &{\rm on}\ \ S_0.\cr}
\leqno(2.3)
$$
The boundary condition $(2.3)_2$ implies that $\psi=\const$ on $S_0$. Hence, 
in view of the definition of $\psi$, we can assume that $\psi=0$ on $S_0$. 
Therefore (2.3) implies the Dirichlet problem
$$\eqal{
&\Delta^{(2)}\psi=b\quad &{\rm in}\ \ \Omega,\cr
&\psi=0\quad &{\rm on}\ \ S_0.\cr}
\leqno(2.4)
$$
Problem (2.4) yields existence of $\psi$ in $H^{s+2}(\Omega)$ and the estimate
$$
\|\psi\|_{H^{s+2}(\Omega)}\le c\|b\|_{H^s(\Omega)},
$$
so (2.2) holds. This concludes the proof.
\kwadrat

\proclaim Lemma 2.2. 
Let $\nabla u\in L_2(\Omega)$ and $u|_{S_0}=0$. Then the following Poincar\'e 
inequality holds
$$
c_p\|u\|_{L_2(\Omega)}\le\|\nabla u\|_{L_2(\Omega)}.
\leqno(2.5)
$$

\noindent
Let us consider the elliptic overdetermined problem
$$\eqal{
&\rot u=b\quad &{\rm in}\ \ D,\cr
&\divv u=0\quad &{\rm in}\ \ D,\cr
&u\cdot\bar n=0\quad &{\rm on}\ \ S,\cr}
\leqno(2.6)
$$
where $D$ is a cylinder and the following compatibility condition holds
$$
\divv b=0.
\leqno(2.7)
$$

\proclaim Lemma 2.3. (see also \cite{21}) 
Let $b\in H^i(D)$, $i=1,2$, and satisfy (2.7). Then there exists a solution to 
(2.6) such that $u\in H^{i+1}(D)$ and
$$
\|u\|_{H^{i+1}(D)}\le c_e\|b\|_{H^i(D)},\quad i=0,1,\quad H^0(D)=L_2(D),
\leqno(2.8)
$$
where $c_e$ does not depend on $u$.

\Proof 
By Lemma 1 from \cite{6}, $(2.6)_{2,3}$ imply existence of a vector $e$ 
such that
$$
u=\rot e,\quad \divv e=0,\quad e_\tau|_S=0.
\leqno(2.9)
$$
The explicit construction of such  vector is presented in \cite{6} and also in 
[23, Sect. 3]. In view of (2.9) problem (2.6) takes the form
$$
-\Delta e=b,\quad e_\tau|_S=0,\quad \divv e|_S=0,
\leqno(2.10)
$$
where $e_\tau=e\cdot\bar\tau$.
The second boundary condition in (2.10) guarantees that $\divv e=0$ in $D$.

\noindent
Recalling the normal and tangent vectors to $S$,
$$\eqal{
&\bar n|_{S_1}={\nabla\varphi\over|\nabla\varphi|},\quad 
\bar\tau_1|_{S_1}={\nabla^\perp\varphi\over|\nabla\varphi|},\quad 
\bar\tau_2|_{S_1}=(0,0,1),\cr
&\bar n|_{S_2}=(0,0,1),\quad \bar\tau_1|_{S_2}=(1,0,0),\quad 
\bar\tau_2|_{S_2}=(0,1,0),\cr}
$$
we express problem (2.10) in the form
$$\eqal{
&-\Delta e=b\quad &{\rm in}\ \ D,\cr
&e_{\tau_1}=0,\ \ e_{\tau_2}=0,\ \ \bar n\cdot\nabla e_n+e_n\divv\bar n=0\quad 
&{\rm on}\ \ S_1,\cr
&e_1=0,\ \ e_2=0,\ \ {\partial\over\partial x_3}e_3=0\quad &{\rm on}\ \ S_2,\cr}
\leqno(2.11)
$$
where $e_n=e\cdot\bar n$. To obtain the last boundary conditions on $S_1$ and 
$S_2$ we formulate dive in the curvilinear coordinates corresponding to 
vectors $\bar n$, $\bar\tau_1$, $\bar\tau_2$ and project it on $S$.

To prove existence of solutions to problem (2.11) we use the idea of 
regularizer (see [14, Ch. 4]).
For this we need a partition of unity and appropriate local estimates. To get 
(2.8) we need the local estimates in $H^2$. Such estimates are easily proved 
in neighborhoods of interior points and points on the smooth part of $S$, so 
points located in a positive distance from the edge $\bar S_1\cap\bar S_2$.

\noindent
Since domain $D$ contains right angles between $S_1$ and $S_2$ we are not able 
to obtain the needed estimates in neighborhoods of points of the edge. From 
$(2.11)_3$ it follows that on $S_2$ we have the Dirichlet and the Neumann 
conditions for the Poisson equation. Therefore, we can reflect the solutions 
of the problems with respect to $S_2$. Then the necessary estimates can be 
easily derived. We have to mention that local estimates near $S_1$ are shown 
after its local flattening.

Summarizing, we have existence of solutions to (2.11) such that $e\in H^2(D)$ 
and
$$
\|e\|_{H^2(D)}\le c\|b\|_{L_2(D)}.
$$
This implies (2.8) for $i=0$. Similarly, we have (2.8) for $i=1$. This 
concludes the proof.
\kwadrat

\noindent
Let $u$ satisfy
$$\eqal{
&\divv u=0\quad &{\rm in}\ \ D,\cr
&u\cdot\bar\tau=0\quad &{\rm on}\ \ S.\cr}
\leqno(2.12)
$$

\proclaim Lemma 2.4. (see [6, Lemma 2.1]) 
For any $u$ satisfying (2.12) the inequality holds
$$
\|u\|_{H^{i+1}(D)}\le c\|\rot u\|_{H^i(D)},\quad i=0,1,
\leqno(2.13)
$$
where $c$ does not depend on $u$.

\noindent
We also need the direct and inverse trace theorems for spaces with mixed norms.

\proclaim Lemma 2.5. (see \cite{5}) 
\item{(i)} Let $u\in W_{p,p_0}^{s,s/2}(\Omega^T)$, $s\in\R_+$, $s>2/p_0$, 
$p,p_0\in(1,\infty)$. Then $u(x,t_0)=u(x,t)|_{t=t_0}$ for $t_0\in[0,T]$ 
belongs to $B_{p,p_0}^{s-2/p_0}(\Omega)$ and
$$
\|u(\cdot,t_0)\|_{B_{p,p_0}^{s-2/p_0}(\Omega)}\le c
\|u\|_{W_{p,p_0}^{s,s/2}(\Omega^T)},
\leqno(2.14)
$$
where $c$ does not depend on $u$.
\item{(ii)} For a given $\tilde u\in B_{p,p_0}^{s-2/p_0}(\Omega)$, $s\in\R_+$, 
$s>2/p_0$, $(p,p_0)\in(1,\infty)$, there exists a function 
$u\in W_{p,p_0}^{s,s/2}(\Omega^T)$, such that $u|_{t=t_0}=\tilde u$ for 
$t_0\in[0,T]$ and
$$
\|u\|_{W_{p,p_0}^{s,s/2}(\Omega^T)}\le c
\|\tilde u\|_{B_{p,p_0}^{s-2/p_0}(\Omega)},
\leqno(2.15)
$$
where $c$ does not depend on $u$.

\proclaim Lemma 2.6. 
For sufficiently regular solutions to (1.7) the following formula is valid
$$
\intop_\Omega w\cdot\nabla w\cdot\Delta wdx=-\intop_{S_0}
(w\cdot\nabla w_1n_2-w\cdot\nabla w_2n_1)\rot^{(2)}wdS_0,
\leqno(2.16)
$$
where $n_1$, $n_2$ are the Cartesian coordinates of the normal vector $\bar n$ 
to $S_0$ and the r.h.s. of (2.16) vanishes by $(1.7)_3$.

\Proof
Let $\rot^{(2)}w=w_{2,x_1}-w_{1,x_2}$ and 
$\tilde\rot\varphi=(\varphi_{,x_2}-\varphi_{,x_1})$. Then
$$\eqal{
\tilde\rot\rot^{(2)}w&=\left(\matrix{
\phantom{-}(\rot^{(2)}w)_{,x_2}\cr -(\rot^{(2)}w)_{,x_1}\cr}\right)=\left(\matrix{
\phantom{-}(w_{2,x_1}-w_{1,x_2})_{,x_2}\cr -(w_{2,x_1}-w_{1,x_2})_{,x_1}\cr}\right)\cr
&=\left(\matrix{
w_{2,x_1x_2}-w_{1,x_2x_2}\cr w_{1,x_2x_1}-w_{2,x_1x_1}\cr}\right)\equiv I.
\cr}
$$
Using the continuity equation $(1.7)_2$ we have
$$
w_{2,x_1x_2}=-w_{1,x_1x_1},\quad w_{1,x_1x_2}=-w_{2,x_2x_2}.
$$
Then $I=-\Delta w$, so
$$
\tilde\rot\rot^{(2)}w=-\Delta w.
\leqno(2.17)
$$
Using (2.17) we have
$$\eqal{
&\intop_\Omega w\cdot\nabla w\cdot\Delta wdx=-\intop_\Omega w\cdot\nabla w
\tilde\rot\rot^{(2)}wdx\cr
&=-\intop_\Omega[w\cdot\nabla w_1(\rot^{(2)}w)_{,x_2}-w\cdot\nabla w_2
(\rot^{(2)}w)_{,x_1}]dx\cr
&=\intop_\Omega[(w\cdot\nabla w_2\rot^{(2)}w)_{,x_1}-
(w\cdot\nabla w_1\rot^{(2)}w)_{,x_2}]dx\cr
&\quad-\intop_\Omega[(w\cdot\nabla w_2)_{,x_1}\rot^{(2)}w-
(w\cdot\nabla w_1)_{,x_2}\rot^{(2)}w]dx\equiv I_1+I_2,\cr}
\leqno(2.18)
$$
where
$$
I_1=\intop_{S_0}(w\cdot\nabla w_2n_1-w\cdot\nabla w_1n_2)\rot^{(2)}wdS_0
$$
and
$$\eqal{
I_2&=-\intop_\Omega(w\cdot\nabla w_{2,x_1}-w\cdot\nabla w_{1,x_2})
\rot^{(2)}wdx\cr
&\quad-\intop_\Omega[(w_{1,x_1}w_{2,x_1}+w_{2,x_1}w_{2,x_2}-w_{1,x_2}w_{1,x_1}-
w_{2,x_2}w_{1,x_2})\rot^{(2)}w]dx\cr
&\equiv I_2^1+I_2^2.\cr}
$$
Continuing,
$$\eqal{
I_2^1&=-\intop_\Omega w\cdot\nabla\rot^{(2)}w\,\rot^{(2)}wdx=-{1\over2}
\intop_\Omega w\cdot\nabla(\rot^{(2)}w)^2dx\cr
&={1\over2}\intop_{S_0}w\cdot\bar n(\rot^{(2)}w)^2dS_0=0\cr}
$$
and
$$
I_2^2=-\intop_\Omega\divv w(\rot^{(2)}w)^2dx=0.
$$
Hence (2.16) holds. This concludes the proof.
\kwadrat

\noindent
Since sometimes is more appropriate to use the slip boundary condtions we 
find a relation between the Navier and the slip boundary conditions.

\proclaim Lemma 2.7. 
Let $\D^{(2)}(w)=\{w_{i,x_j}+w_{j,x_i}\}_{i,j=1,2}$, 
$\rot^{(2)}w=w_{2,x_1}-w_{1,x_2}$. Let $\bar n=(n_1,n_2)$ be the normal unit 
outward vector to $S_0$ and $\bar\tau=(-n_2,n_1)$ be the tangent vector. Then
$$
\bar n\cdot\D^{(2)}(w)\cdot\bar\tau=\rot^{(2)}w+2w_{n,\tau}-2w_in_{i,\tau},
\leqno(2.19)
$$
where $w_n=w\cdot\bar n$, $w_{,\tau}=\bar\tau\cdot\nabla w$.

\Proof 
$$\eqal{
&\bar n\cdot\D^{(2)}(w)\cdot\bar\tau=n_i(w_{i,x_j}+w_{j,x_i})\tau_j=n_i
(-w_{i,x_j}+w_{j,x_i})\tau_j\cr
&\quad+2n_iw_{i,x_j}\tau_j=n_1(-w_{1,x_2}+w_{2,x_1})\tau_2+n_2
(-w_{2,x_1}+w_{1,x_2})\tau_1\cr
&\quad+2n_iw_{i,x_j}\tau_j=(n_1\tau_2-n_2\tau_1)\rot^{(2)}w+2w_{n,\tau}-
2w_in_{i,\tau}\cr
&=\rot^{(2)}w+2w_{n,\tau}-2w_in_{i,\tau}.\cr}
$$
This concludes the proof.

Let us consider the Stokes problem
$$\eqal{
&w_t+\nu\tilde\rot\rot^{(2)}w+\nabla\eta=f\quad &{\rm in}\ \ \Omega\times(0,T),
\cr
&\divv w=0\quad &{\rm in}\ \ \Omega\times(0,T),\cr
&w\cdot\bar n=0,\ \ \rot^{(2)}w=0\quad &{\rm on}\ \ S_0\times(0,T),\cr
&w|_{t=0}=w(0)\quad &{\rm in}\ \ \Omega.\cr}
\leqno(2.20)
$$
The theory developed in \cite{12, 13, 19, 20} implies

\proclaim Lemma 2.8. 
Let $f\in L_{\sigma_1,\sigma_2}(\Omega\times(0,T))$, 
$w(0)\in B_{\sigma_1,\sigma_2}^{2-2/\sigma_2}(\Omega)$, 
$\sigma_1,\sigma_2\in(1,\infty)$. Then there exists a solution to problem 
(2.20) such that $w\in W_{\sigma_1,\sigma_2}^{2,1}(\Omega\times(0,T))$, 
$\nabla\eta\in L_{\sigma_1,\sigma_2}(\Omega\times(0,T))$ and
$$\eqal{
&\|w\|_{W_{\sigma_1,\sigma_2}^{2,1}(\Omega\times(0,T))}+
\|\nabla\eta\|_{L_{\sigma_1,\sigma_2}(\Omega\times(0,T))}\cr
&\le c(\|f\|_{L_{\sigma_1,\sigma_2}(\Omega\times(0,T))}+
\|w(0)\|_{B_{\sigma_1,\sigma_2}^{2-2/\sigma_2}(\Omega)}),\cr}
\leqno(2.21)
$$
where $c$ may depend on $T$.

\section{3. Two-dimensional solutions}

First we have

\proclaim Lemma 3.1. 
Assume that $\bar h\in L_{2,loc}(\R_+;L_2(\Omega))$, $w(0)\in L_2(\Omega)$. 
Assume that $T>0$ is given. Denote 
$A_1^2=\sup_{k\in\N_0}{1\over\nu c_1}\intop_{kT}^{(k+1)T}
\|\bar h(t)\|_{L_2(\Omega)}^2dt<\infty$, 
$A_2^2={A_1^2\over1-e^{-\nu c_1T}}+\|w(0)\|_{L_2(\Omega)}^2<\infty$, 
where $c_1$ appearing in (3.5) follows from Lemma 2.1 (see (2.2)). 
Then for solutions to (1.7) we have
$$
\|w(kT)\|_{L_2(\Omega)}^2\le A_2^2
\leqno(3.1)
$$
and
$$
\|w(t)\|_{L_2(\Omega)}^2+\nu c_1\intop_{kT}^t\|w(t')\|_{H^1(\Omega)}^2dt'\le
A_1^2+A_2^2\equiv A_3^2,
\leqno(3.2)
$$
where $t\in(kT,(k+1)T]$ and $k\in\N_0$.

\Proof 
Multiplying $(1.7)_1$ by $w$ and integrating over $\Omega$ yields
$$
{1\over2}{d\over dt}\|w\|_{L_2(\Omega)}^2+\nu\intop_\Omega\tilde\rot\rot^{(2)}w
\cdot wdx=\intop_\Omega\bar h\cdot wdx,
\leqno(3.3)
$$
where the first of boundary conditions $(1.7)_3$ is used.

\noindent
The second term on the l.h.s. of (3.3) equals
$$\eqal{
&\intop_\Omega[(\rot^{(2)}w)_{,x_2}w_1-(\rot^{(2)}w)_{,x_1}w_2]dx\cr
&=\intop_\Omega[(\rot^{(2)}ww_1)_{,x_2}+(-\rot^{(2)}ww_2)_{,x_1}]dx+ \intop_\Omega|\rot^{(2)}w|^2dx.\cr}
\leqno(3.4)
$$
Applying the Green formula, the first term on the r.h.s. of (3.4) is equal to
$$
\intop_{S_0}\rot^{(2)}w(w_1n_2-w_2n_1)dS_0=0,
$$
where the second condition from $(1.7)_3$ is utilized.

\noindent
Employing (3.4), Lemma 2.1 and the H\"older and the Young inequalities to the 
r.h.s. of (3.3) we obtain from (3.3) the inequality
$$
{d\over dt}\|w\|_{L_2(\Omega)}^2+\nu c_1\|w\|_{H^1(\Omega)}^2\le{1\over\nu c_1}
\|\bar h\|_{L_2(\Omega)}^2.
\leqno(3.5)
$$
Expressing (3.5) in the form
$$
{d\over dt}\|w\|_{L_2(\Omega)}^2+\nu c_1\|w\|_{L_2(\Omega)}^2\le{1\over\nu c_1}
\|\bar h\|_{L_2(\Omega)}^2
\leqno(3.6)
$$
we integrate it with respect to time from $t=kT$ to $t\in(kT,(k+1)T]$, 
$k\in\N_0$, to derive
$$\eqal{
\|w(t)\|_{L_2(\Omega)}^2&\le{1\over\nu c_1}\intop_{kT}^t
\|\bar h(t')\|_{L_2(\Omega)}^2dt'\cr
&\quad+\|w(kT)\|_{L_2(\Omega)}^2\exp(-\nu c_1(t-kT)).\cr}
\leqno(3.7)
$$
Setting $t=(k+1)T$ inequality (3.7) implies
$$
\|w((k+1)T)\|_{L_2(\Omega)}^2\le{1\over\nu c_1}\intop_{kT}^t
\|\bar h(t)\|_{L_2(\Omega)}^2dt+\|w(kT)\|_{L_2(\Omega)}^2e^{-\nu c_1T}.
\leqno(3.8)
$$
By iteration we get
$$
\|w(kT)\|_{L_2(\Omega)}^2\le{A_1^2\over1-e^{-\nu c_1T}}+
\|w(0)\|_{L_2(\Omega)}^2e^{-\nu c_1kT}\le A_2^2,
\leqno(3.9)
$$
so (3.1) holds.
Integrating (3.5) with respect to time from $t=kT$ to $t\in(kT,(k+1)T]$ and 
using (3.1) yields (3.2). This concludes the proof.
\kwadrat

\proclaim Lemma 3.2. 
Let the assumptions of Lemma 3.1 hold. Let $w(0)\in H^1(\Omega)$. Let 
$A_4^2={c_1A_1^2\over1-e^{-c_p\nu T}}+\|\rot^{(2)}(0)\|_{L_2(\Omega)}^2$. Then
$$
\|\rot^{(2)}w(kT)\|_{L_2(\Omega)}^2\le A_4^2
\leqno(3.10)
$$
and
$$\eqal{
&\|\rot^{(2)}w(t)\|_{L_2(\Omega)}^2+\nu c_p\intop_{kT}^t
\|\rot^{(2)}w(t')\|_{H^1(\Omega)}^2dt'\cr
&\le{1\over\nu}\intop_{kT}^t\|\bar h(t')\|_{L_2(\Omega)}^2dt'+
\|\rot^{(2)}w(kT)\|_{L_2(\Omega)}^2\le c_1A_1^2+A_4^2\equiv A_5^2.\cr}
\leqno(3.11)
$$

\Proof 
Multiplying $(1.7)_1$ by $-\Delta w$ and integrating the result over $\Omega$ 
yields
$$\eqal{
&-\intop_\Omega w_{,t}\cdot\Delta wdx+\nu\intop_\Omega|\Delta w|^2dx=
\intop_\Omega\nabla\bar\eta\cdot\Delta wdx+\intop_\Omega w\cdot\nabla w\cdot
\Delta wdx\cr
&\quad-\intop_\Omega\bar h\cdot\Delta wdx.\cr}
\leqno(3.12)
$$
Using that $\Delta w=-\tilde\rot\rot^{(2)}w$ the first term on the l.h.s. of 
(3.12) equals
$$\eqal{
&\intop_\Omega w_t\tilde\rot\rot^{(2)}wdx=\intop_\Omega
[w_{1,t}\partial_{x_2}\rot^{(2)}w-w_{2,t}\partial_{x_1}\rot^{(2)}w]dx\cr
&=\intop_\Omega[(-w_{2,t}\rot^{(2)}w)_{,x_1}+
(w_{1,t}\rot^{(2)}w)_{,x_2}]dx\cr
&\quad+\intop_\Omega[w_{2,x_1t}\rot^{(2)}w-w_{1,x_2t}\rot^{(2)}w]dx=
\intop_\Omega\rot^{(2)}w\partial_t\rot^{(2)}wdx\cr
&={1\over2}{d\over dt}\intop_\Omega|\rot^{(2)}w|^2dx,\cr}
$$
where we used the boundary conditions $(1.7)_3$. In view of the boundary 
conditions $(1.7)_3$ also and Lemma 2.6 the second term on the r.h.s. of 
(3.12) vanishes. The first term on the r.h.s. of (3.12) takes the form
$$\eqal{
&-\intop_\Omega[\partial_{x_1}\bar\eta(\rot^{(2)}w)_{,x_2}-
\partial_{x_2}\bar\eta(\rot^{(2)}w)_{,x_1}]dx\cr
&=\intop_\Omega[\partial_{x_2}(\partial_{x_1}\bar\eta\rot^{(2)}w)+
\partial_{x_1}(-\partial_{x_2}\bar\eta\rot^{(2)}w)]dx\cr}
$$
which also vanishes in view of boundary conditions $(1.7)_3$.

\noindent
Using the above calculations in (3.12) yields
$$
{1\over2}{d\over dt}\|\rot^{(2)}w\|_{L_2(\Omega)}^2+\nu
\|\nabla\rot^{(2)}w\|_{L_2(\Omega)}^2=\intop_\Omega\bar h\cdot\Delta wdx
\leqno(3.13)
$$
Applying the H\"older and the Young inequalitites to the r.h.s. of (3.13) 
and using that $|\Delta w|=|\nabla\rot^{(2)}w|$ we obtain
$$
{d\over dt}\|\rot^{(2)}w\|_{L_2(\Omega)}^2+\nu
\|\nabla\rot^{(2)}w\|_{L_2(\Omega)}^2\le{1\over\nu}\|\bar h\|_{L_2(\Omega)}^2
\leqno(3.14)
$$
Since $\rot^{(2)}w|_{S_0}=0$ we can apply the Poincar\'e inequality (see 
(2.5)) to (3.14). Hence, we get
$$
{d\over dt}\|\rot^{(2)}w\|_{L_2(\Omega)}^2+c_p\nu
\|\rot^{(2)}w\|_{L_2(\Omega)}^2\le{1\over\nu}\|\bar h\|_{L_2(\Omega)}^2,
\leqno(3.15)
$$
where $c_p$ is the constant from the Poincar\'e inequality (2.5). Integrating 
(3.15) with respect to time from $t=kT$ to $t\in(kT,(k+1)T]$, $k\in\N_0$, 
we derive
$$\eqal{
&\|\rot^{(2)}w(t)\|_{L_2(\Omega)}^2\le{1\over\nu}e^{-c_p\nu T}\intop_{kT}^t
\|\bar h(t')\|_{L_2(\Omega)}^2e^{c_p\nu t'}dt'\cr
&\quad+\|\rot^{(2)}w(kT)\|_{L_2(\Omega)}^2\exp(-c_p\nu(t-kT)).\cr}
\leqno(3.16)
$$
Setting $t=(k+1)T$ in (3.16) yields
$$\eqal{
&\|\rot^{(2)}w((k+1)T)\|_{L_2(\Omega)}^2\le{1\over\nu}\intop_{kT}^{(k+1)T}
\|\bar h(t)\|_{L_2(\Omega)}^2dt\cr
&\quad+\|\rot^{(2)}w(kT)\|_{L_2(\Omega)}^2\exp(-c_p\nu T).\cr}
\leqno(3.17)
$$
By iteration we have
$$
\|\rot^{(2)}w(kT)\|_{L_2(\Omega)}^2\le{c_1A_1^2\over1-e^{-c_p\nu T}}+
\|\rot^{(2)}w(0)\|e^{-c_p\nu kT}\le A_4^2
\leqno(3.18)
$$
Hence (3.10) holds.
Integrating (3.14) with respect to time from $t=kT$ to $t\in(kT,(k+1)T]$, 
$k\in\N_0$, using (3.10) and the Poincar\'e inequality (2.5) we derive (3.11). 
This concludes the proof.
\kwadrat

\Remark{3.3.} 
In view of Lemma 2.1 inequalities (3.10) and (3.11) can be expressed in the 
form
$$
\|w(kT)\|_{H^1(\Omega)}^2\le c_2A_4^2
\leqno(3.19)
$$
and
$$
\|w(t)\|_{H^1(\Omega)}^2+\nu c_p\intop_{kT}^t\|w(t')\|_{H^2(\Omega)}^2dt'\le
c_2A_5^2,
\leqno(3.20)
$$
where $c_2$ depends on the constant $c$ from (2.2). Proof of existence is 
standard.

From Lemmas 3.1, 3.2 and Remark 3.3, Theorem 1 follows.

\noindent
We prove that $w\in C([kT,(k+1)T];W_\sigma^1(\Omega))$, $\sigma>3$, $k\in\N_0$. Hence, we want to show that
$$
\|w(t)\|_{W_\sigma^1(\Omega)}\le A_6,
\leqno(3.21)
$$
where $A_6$ does not depend on time.

The above increasing of regularity is made in \cite{22} by the applying the 
energy method. This needs much more regularity of data than it is necessary 
to show (3.21). Moreover, it implies a stronger relation between dissipation 
and the external force than it is presented in (4.2). Therefore, we follow 
the regularity increasing technique used in \cite{24}.

\noindent
In this case we have only restriction (4.2). The above mentioned method from 
\cite{24} is possible because Lemma 3.2 and Remark 3.3 imply that 
$w\cdot\nabla w\in L_2(kT,(k+1)T;L_\sigma(\Omega))$, $\sigma\in(3,\infty)$ 
and $k\in\N_0$.

\proclaim Lemma 3.4. 
Assume that $w(0)\in B_{\sigma,2}^2(\Omega)$, 
$\bar h\in L_2(kT,(k+1)T;L_\sigma(\Omega))$ $k\in\N_0$, $\sigma>3$. Then 
$w\in C(\R_+;W_\sigma^1(\Omega))$, $\sigma>3$ and (3.21) holds with constant 
$A_6$ depending on $\|w(0)\|_{B_{\sigma,2}^2(\Omega)}$ and 
$\sup_k\intop_{kT}^{(k+1)T}\|\bar h(t)\|_{L_\sigma(\Omega)}dt$.

\Proof 
Remark 3.3 implies that $w\cdot\nabla w\in L_2(kT,(k+1)T;L_\sigma(\Omega))$, 
$\sigma>3$. In view of the assumptions of the lemma, the theory developed in 
\cite{12, 13, 19, 20} implies existence of solutions to 
problem (1.7) such that $w\in W_{\sigma,2}^{2,1}(\Omega\times(kT,(k+1)T))$, 
$\nabla\bar\eta\in L_{\sigma,2}(\Omega\times(kT,(k+1)T))$ and
$$\eqal{
&\|w\|_{W_{\sigma,2}^{2,1}(\Omega\times(kT,(k+1)T))}\le c(A_5^2+
\|\bar h\|_{L_2(kT,(k+1)T;L_\sigma(\Omega))}\cr
&\quad+\|w(kT)\|_{B_{\sigma,2}^2(\Omega)}),\cr}
\leqno(3.22)
$$
where $c$ may depend on $T$. Inequality (3.22) implies (3.21) if we know that 
$\|w(kT)\|_{B_{\sigma,2}^1(\Omega)}$ is bounded by a constant independent of 
$k$. Hence for $k=0$, (3.22) implies (3.21). For $k=1$ and Lemma 2.5 we 
calculate $\|w(T)\|_{B_{\sigma,2}^1(\Omega)}$ using (3.22) for $k=0$.
Then Lemma 2.8 implies existence of solutions to (1.7) such that
$$
\|w\|_{W_{\sigma,2}^{2,1}(\Omega\times(T,2T))}\le c,
$$
where $c$ depends on $\|w(0)\|_{B_{\sigma,2}^1(\Omega)}$ and 
$\|\bar h\|_{L_2(kT,(k+1)T;L_\sigma(\Omega))}$ for $k=0,1$.

To eliminate dependence on $\|w(kT)\|_{B_{\sigma,2}^1(\Omega)}$ in r.h.s. of 
(3.22) we use a smooth cut-off function $\zeta_k=\zeta_k(t)$ such that 
$\zeta_k(t)=0$ for $t\in[kT,kT+\delta/2]$ and $\zeta_k(t)=1$ for 
$t\in[kT+\delta,(k+1)T]$, where $\delta<T$. Introducing the quantities
$$
w_k=w\zeta_k,\quad \bar\eta_k=\bar\eta\zeta_k,\quad \bar h_k=\bar h\zeta_k
$$
we see that problem (1.7) takes the form
$$\eqal{
&w_{k,t}+w\cdot\nabla w_k+\nu\tilde\rot\rot^{(2)}w_k+\nabla\bar\eta_k=
w\dot\zeta_k+\bar h_k,\cr
&\divv w_k=0,\cr
&w_k\cdot\bar n=0,\ \ \rot^{(2)}w_k=0\quad {\rm on}\ \ S_0,\cr
&w_k|_{t=kT}=0.\cr}
\leqno(3.23)
$$
In view of Lemmas 2.8 and 3.2 we obtain for solutions to (3.23) the estimate
$$\eqal{
&\|w_k\|_{W_{\sigma,2}^{2,1}(\Omega\times(kT+\delta/2,(k+1)T))}\cr 
&\le c(A_5+A_5^2+\|\bar h_k\|_{L_{\sigma,2}(\Omega\times(kT+\delta/2,(k+1)T))})\equiv cA_7,\cr}
\leqno(3.24)
$$
where in view of the estimate for $\bar h$ (see assumptions of Lemma 3.2) we 
see that $A_7$ does not depend on $k$. Then Lemma 2.5 implies
$$
\|w((k+1)T)\|_{B_{\sigma,2}^1(\Omega)}\le cA_7.
\leqno(3.25)
$$
Applying Lemma 2.8 and using (3.25) we obtain
$$\eqal{
&\|w\|_{W_{\sigma,2}^{2,1}(\Omega\times((k+1)T,(k+1)T+\delta/2))}\cr
&\le c(A_5+A_5^2+\|\bar h_k\|_{L_{\sigma,2}(\Omega\times(kT,(k+2)T))})\le cA_7
\cr}
\leqno(3.26)
$$
In view of (3.24) and (3.26) we prove (3.21) with constant $A_6$ independent 
of $k$. This concludes the proof.

\Remark{3.5}
Remark 3.3 implies that 
$\|w\cdot\nabla w\|_{L_2(\Omega\times(kT,(k+1)T))}\le cA_5^2$, so the regularizer technique (see \cite[Ch. 4]{14}) gives existence of solutions to problem (1.5) such that $w\in W_2^{2,1}(\Omega\times(kT,(k+1)T))$, $\nabla\eta\in L_2(\Omega\times(kT,(k+1)T))$, $k\in\N_0$ and
$$
\|w\|_{W_2^{2,1}(\Omega\times(kT,(k+1)T))}+
\|\nabla\eta\|_{L_2(\Omega(kT,(k+1)T))}\le c(A_5+A_5^2).
$$

\section{4. Stability}

In this Section we prove stability of two-dimensional solutions. For this 
purpose we examine problem (1.9). First we show the $L_2$-stability.

\proclaim Lemma 4.1. 
Let the assumptions of Lemma 3.2 be satisfied. Assume that
$$\eqal{
&B_1^2=\sup_k\intop_{kT}^{(k+1)T}\|g(t)\|_{L_2(D)}^2dt,\quad
B_2^2={c\over\nu}\exp\bigg({c\over\nu}A_5^2\bigg)B_1^2,\cr
&B_3^2={B_2^2\over1-\exp\big(-{\nu\over2}T\big)}+\|u(0)\|_{L_2(D)}^2.\cr}
\leqno(4.1)
$$
$$
-{\nu\over2}T+{c\over\nu}A_5^2\le0.
\leqno(4.2)
$$
Then the following estimates for solutions to (1.9) hold
$$
\|u(kT)\|_{L_2(D)}^2\le B_3^2
\leqno(4.3)
$$
$$
\|u(t)\|_{L_2(D)}^2\le\exp\bigg({c\over\nu}A_5^2\bigg){c\over\nu}B_1^2+B_3^2
\equiv B_4^2.
\leqno(4.4)
$$

\Proof 
Multiplying $(1.9)_1$ by $u$, integrating over $D$ and using the boundary 
conditions yield
$$
{1\over2}{d\over dt}\|u\|_{L_2(D)}^2+\nu\|\rot u\|_{L_2(D)}^2=-\intop_D
u\cdot\nabla w\cdot udx+\intop_D g\cdot udx
\leqno(4.5)
$$
Applying the H\"older and the Young inequalities to the r.h.s. of (4.5) and 
using Lemma 2.1 we obtain
$$\eqal{
&{1\over2}{d\over dt}\|u\|_{L_2(D)}^2+{\nu\over2}\|\rot u\|_{L_2(D)}^2\le
{c_e\over\nu}\|\nabla w\|_{L_3(D)}^2\|u\|_{L_2(D)}^2\cr
&\quad+{1\over\nu}\|g\|_{L_2(D)}^2,\cr}
\leqno(4.6)
$$
where $c_e$ appears in (2.8). Applying again Lemma 2.1 we have
$$\eqal{
&{d\over dt}\|u\|_{L_2(D)}^2+\nu\|u\|_{L_2(D)}^2\le{c\over\nu}
\|\nabla w\|_{L_3(D)}^2\|u\|_{L_2(D)}^2\cr
&\quad+{c\over\nu}\|g\|_{L_2(D)}^2.\cr}
\leqno(4.7)
$$
Inequality (4.7) implies
$$\eqal{
&{d\over dt}\bigg[\|u\|_{L_2(D)}^2\exp\bigg(\nu t-{c\over\nu}\intop_{kT}^t
\|\nabla w(t')\|_{L_3(D)}^2dt'\bigg)\bigg]\cr
&\le{c\over\nu}\|g\|_{L_2(D)}^2\exp\bigg(\nu t-{c\over\nu}\intop_{kT}^t
\|\nabla w(t')\|_{L_3(D)}^2dt'\bigg).\cr}
\leqno(4.8)
$$
Integrating (4.8) with respect to time from $t=kT$ to $t\in(kT,(k+1)T]$ gives
$$\eqal{
&\|u(t)\|_{L_2(D)}^2\le\exp\bigg(-\nu t+{c\over\nu}\intop_{kT}^t
\|\nabla w(t')\|_{L_3(D)}^2dt'\bigg)\cdot\cr
&\quad\cdot{c\over\nu}\intop_{kT}^t\|g(t')\|_{L_2(D)}^2\exp(\nu t')dt'\cr
&\quad+\exp\bigg(-\nu(t-kT)+{c\over\nu}\intop_{kT}^t
\|\nabla w(t')\|_{L_3(D)}^2dt'\bigg)\|u(kT)\|_{L_2(D)}^2.\cr}
\leqno(4.9)
$$
Setting $t=(k+1)T$ and employing (3.11) we have
$$\eqal{
&\|u((k+1)T)\|_{L_2(D)}^2\le\exp\bigg({c\over\nu}A_5^2\bigg){c\over\nu}
\intop_{kT}^{(k+1)T}\|g(t')\|_{L_2(D)}^2dt'\cr
&\quad+\exp\bigg(-\nu T+{c\over\nu}A_5^2\bigg)\|u(kT)\|_{L_2(D)}^2.\cr}
\leqno(4.10)
$$
In view of assumptions (4.1) inequality (4.10) takes the form
$$
\|u((k+1)T)\|_{L_2(D)}^2\le B_2^2+\exp\bigg(-{\nu\over2}T\bigg)
\|u(kT)\|_{L_2(D)}^2.
\leqno(4.11)
$$
By iteration we have
$$
\|u(kT)\|_{L_2(D)}^2\le{B_2^2\over1-\exp\big(-{\nu\over2}T\big)}+\exp
\bigg(-{\nu\over2}kT\bigg)\|u(0)\|_{L_2(D)}^2\le B_3^2,
\leqno(4.12)
$$
so (4.3) holds. Applying (4.12), (4.2) and (4.3) in (4.9) yields (4.4). 
This concludes the proof.

Finally, we prove stability omitting the strong restriction (4.2).

\proclaim Lemma 4.2. 
Let $w\in L_2(kT,(k+1)T;W_{2^+}^1(\Omega))$, $g\in C(kT,(k+1)T;\break L_2(D))$, 
$k\in\N_0$. Let $c_*$ be a constant such that 
$c_*\in(0,\nu]$. Let $\gamma_*$ be so small that 
${c_*\over2}\le\nu-{c_0\over\nu^3}\gamma_*^4$, where $c_0$ is the constant 
from (4.20). Let $\gamma\le\gamma_*$. Let
$$\eqal{
&\|u(0)\|_{H^1(\Omega)}^2\le\gamma\cr
&G^2(t)={c\over\nu}\|g(t)\|_{L_2(D)}^2\le{c_*\over4}\gamma.\cr}
\leqno(4.13)
$$
Then
$$
\|u(t)\|_{H^1(\Omega)}^2\le\gamma\quad {\sl for}\ \ t\in\R_+.
\leqno(4.14)
$$

\Proof 
Multiplying (1.9) by $\rot^2u$, integrating over $D$ and by parts yields
$$\eqal{
&{1\over2}{d\over dt}\|\rot u\|_{L_2}^2+\nu\|\rot^2u\|_{L_2}^2\le-\intop_D
u\cdot\nabla u\cdot\rot^2udx\cr
&\quad-\intop_Dw\cdot\nabla u\cdot\rot^2udx-\intop_Du\cdot\nabla w\cdot\rot^2udx+\intop_Dg\rot^2udx.\cr}
\leqno(4.15)
$$
The first term on the r.h.s. of (4.15) equals
$$\eqal{
&-\intop_\Omega\rot(u\cdot\nabla u)\rot udx\cr 
&=-\intop_\Omega u\cdot\nabla\rot u
\rot udx-\intop_\Omega\varepsilon_{kij}u_{l,x_j}\partial_{x_l}u_i
(\rot u)_kdx\equiv I_1,\cr}
$$
where $\varepsilon_{kij}$ is the antisymmetric Ricci tensor and summation is 
performed over all repeated indices.
Since the first term in $I_1$ vanishes in view of the boundary conditions, 
we have
$$
|I_1|\le c\|u_x\|_{L_3(D)}^3.
$$
Applying the H\"older and the Young inequalities to the other terms on the 
r.h.s. of (4.15), we obtain
$$\eqal{
&{d\over dt}\|\rot u\|_{L_2(D)}^2+\nu\|\rot^2u\|_{L_2(D)}^2\le c
\|u_x\|_{L_3(D)}^3+{c\over\nu}\|w\cdot\nabla u\|_{L_2(D)}^2\cr
&\quad+{c\over\nu}\|u\cdot\nabla w\|_{L_2(D)}^2+{c\over\nu}
\|g\|_{L_2(D)}^2.\cr}
\leqno(4.16)
$$
In view of the interpolation inequality
$$
\|u_x\|_{L_3(D)}\le c\|u_x\|_{H^1(D)}^{1/2}\|u_x\|_{L_2(D)}^{1/2}
$$
the first term on the r.h.s. of (4.16) is bounded by
$$
c\|u_x\|_{H^1(D)}^{3/2}\|u_x\|_{L_2(D)}^{3/2}\le{\varepsilon_1^{4/3}\over4/3}
\|u_x\|_{H^1(D)}^2+{c\over4\varepsilon_1^4}\|u_x\|_{L_2(D)}^6\equiv I_2.
$$
Using Lemma 2.3 we have
$$
\|u_x\|_{L_2(D)}\le\|u\|_{H^1(D)}\le c\|\rot u\|_{L_2(D)}
$$
and Lemma 2.4 gives
$$
\|u_x\|_{H^1(D)}\le c\|\rot u\|_{H^1(D)}\le c\|\rot^2u\|_{L_2(D)}.
\leqno(4.17)
$$
Setting ${\varepsilon_1^{4/3}\over4/3}={\nu\over 4c}$ we obtain that
$$
I_2\le{\nu\over4}\|\rot^2u\|_{L_2(D)}^2+{c\over\nu^3}\|\rot u\|_{L_2(D)}^6.
$$
The second term on the r.h.s. (4.16) is bounded by
$$
{c\over\nu}\|w\|_{L_\infty(D)}^2\|\nabla u\|_{L_2(D)}^2\le{c\over\nu}
\|w\|_{W_{2^+}^1(\Omega)}^2\|\nabla u\|_{L_2(D)}^2,
$$
where $2^+>2$ but is very close to 2.

\noindent
Finally, the third term on the r.h.s. of (4.16) is estimated by
$$
{c\over\nu}\|u\|_{L_6(D)}^2\|\nabla w\|_{L_3(\Omega)}^2\le{c\over\nu}
\|u\|_{H^1(D)}^2\|\nabla w\|_{L_3(\Omega)}^2.
$$
Employing the above estimates in (4.16) yields the inequality
$$\eqal{
&{d\over dt}\|\rot u\|_{L_2(D)}^2+\nu\|\rot^2u\|_{L_2(D)}^2\le{c\over\nu^3}
\|\rot u\|_{L_2(D)}^6\cr
&\quad+{c\over\nu}\|w\|_{W_{2^+}^1(\Omega)}^2\|u\|_{H^1(D)}^2+
{c\over\nu}\|g\|_{L_2(D)}^2.\cr}
\leqno(4.18)
$$
Let us introduce the quantities
$$\eqal{
&X(t)=\|\rot u(t)\|_{L_2(D)},\quad Y(t)=\|\rot^2u(t)\|_{L_2(D)},\cr
&G^2(t)={c\over\nu}\|g(t)\|_{L_2(D)}^2,\quad
A^2(t)={c\over\nu}\|w\|_{W_{2^+}^1(\Omega)}^2\cr}
\leqno(4.19)
$$
To prove the lemma we need to know that $G(t)\in C([kT,(k+1)T])$, 
$A^2\in L_1(kT,(k+1)T)$, $k\in\N_0$. 
By the assumptions of the lemma we know that 
$\|g(t)\|_{L_2(D)}\in C([kT,(k+1)T])$. Similarly, Lemma 4.1 implies that 
$\|u(t)\|_{L_2(D)}\in C([kT,(k+1)T])$. Moreover, Lemma 3.4 shows that 
$w\in C([kT,(k+1)T];W_{2^+}^1(\Omega))$ for all $k\in\N_0$.

The aim of this lemma is to show that $\|u(t)\|_{H^1(D)}$ is as small as 
$\|u(0)\|_{H^1(D)}$ in each time interval $[kT,(k+1)T]$, $k\in\N_0$, 
separately. Therefore, in view of notation (4.19) and with the above shown 
properties of $G(t)$, we express (4.18) in the short form
$$
{d\over dt}X^2+\nu Y^2\le{c_0\over\nu^3}X^6+A^2X^2+G^2.
\leqno(4.20)
$$
Since $X\le Y$ we have
$$
{d\over dt}X^2\le-X^2\bigg(\nu-{c_0\over\nu^3}X^4\bigg)+A^2X^2+G^2.
\leqno(4.21)
$$
Let $\gamma\in(0,\gamma_*]$, where $\gamma_*$ is so small that
$$
\nu-{c_0\over\nu^3}\gamma_*^4\ge{c_*\over2},\quad 0<c_*\le\nu.
\leqno(4.22)
$$
Since the coefficients of (4.21) depend on the two-dimensional solutions 
determined step by step in time, we consider (4.21) in the interval 
$[kT,(k+1)T]$, $k\in\N_0$, with the assumptions
$$
X^2(kT)\le\gamma,\quad G^2(t)\le c_*{\gamma\over4},\quad t\in[kT,(k+1)T].
\leqno(4.23)
$$
Let us introduce the quantity
$$
Z^2(t)=\exp\bigg(-\intop_{kT}^tA^2(t')dt'\bigg)X^2(t),\quad t\in[kT,(k+1)T].
$$
Then (4.21) takes the form
$$
{d\over dt}Z^2\le-\bigg(\nu-{c_0\over\nu^3}X^4\bigg)Z^2+\bar G^2,
\leqno(4.24)
$$
where $\bar G^2=G^2\exp\big(-\intop_{kT}^tA^2(t')dt'\big)$.

\noindent
Suppose that
$$\eqal{
t_*&=\inf\{t\in(kT,(k+1)T]:\ X^2(t)>\gamma\}\cr
&=\inf\bigg\{t\in(kT,(k+1)T]:\ Z^2(t)>\gamma\exp
\bigg(-\intop_{kT}^tA^2(t')dt'\bigg)\bigg\}>kT.\cr}
$$
By (4.22) for $t\in(0,t_*]$ inequality (4.24) takes the form
$$
{d\over dt}Z^2\le-{c_*\over2}Z^2+\bar G^2(t)
\leqno(4.25)
$$
Clearly, we have
$$
Z^2(t_*)=\gamma\exp\bigg(-\intop_{kT}^{t_*}A^2(t')dt'\bigg)\quad {\rm and}\quad
Z^2(t)>\gamma\exp\bigg(-\intop_{kT}^tA^2(t')dt'\bigg)
\leqno(4.26)
$$
for $t>t_*$. But (4.23) and (4.25) yield
$$
{d\over dt}Z^2|_{t=t_*}\le c_*\bigg(-{\gamma\over2}+{\gamma\over4}\bigg)\exp
\bigg(-\intop_{kT}^{t_*}A^2(t')dt'\bigg)<0
$$
so it contradicts to (4.26). Hence $X^2((k+1)T)<\gamma$. Then induction proves 
the lemma.
\kwadrat

\Remark{4.3}
Let the assumptions of Remark 3.3 and Lemma 4.2 hold. Then the regularizer technique (see \cite[Ch 4]{14}) gives existence of solutions to problem (1.9) such that $u\in W_2^{2,1}(D\times(kT,(k+1)T))$, 
$\nabla q\in L_2(D\times(kT,(k+1)T))$, $k\in\N_0$, and
$$
\|u\|_{W_2^{2,1}(D\times(kT,(k+1)T))}+\|\nabla q\|_{L_2(D\times(kT,(k+1)T))}\le
c\gamma[A_5+\gamma+1].
$$

\section{References}

\item{1.} Auscher P., Dubois S. and Tchamitchian P.: On the stability of 
global solutions to Navier-Stokes equations in the space, Journal de 
Math\'ematiques Pures et Appliques 83 (2004), 673--697.

\item{2.} Bardos C., Lopes Filho M. C., Niu D., Nussenzveig Lopes H. J. and 
Titi E. S.: Stability of two-dimensional viscous incompressible flows under 
three-dimesional perturbations and inviscid symmetry breaking, SIAM J. Math. 
Anal. 45 (2013), 1871--1885.

\item{3.} Beir\~ao da Veiga H. and Secchi P.: $Lp$-stability for the strong 
solutions of the Navier-Stokes equations in the whole space, Arch. Ration. 
Mech. Anal. 98 (1987), 65--69.

\item{4.} Besov, O. V.; Il'in, V. P.; Nikol'skii, S. M.: Integral 
representations of functions and imbedding theorems, Nauka, Moscow 1975 
(in Russian).

\item{5.} Bugrov, Y. S.: Function spaces with mixed norm, Math. USSR -- 
Izv. 5 (1971), 1145--1167 (in Russian).

\item{6.} Bykhovsky, E. B.: Solvability of mixed problem for the Maxwell 
equations for ideal conductive boundary, Vest. Lenin. Univ. Ser. Mat. Mekh. 
Astron. 13 (1957), 50--66 (in Russian).

\item{7.} Chemin J. I, and Gallagher I.: Wellposedness and stability results 
for the Navier-Stokes equations in $R^3$, Ann. I. H. Poincar\'e, Analyse Non 
Lin\'eaire, 26 (2009), 599--624.

\item{8.} Gallagher I.: The tridimensional Navier-Stokes equations with 
almost bidimensional data: stability, uniqueness and life span, Internat. Mat. 
Res. Notices 18 (1997), 919--935.

\item{9.} Gallagher I., Iftimie D. and Planchon F.: Asymptotics and stability 
for global solutions to the Navier-Stokes equations, Annales de I'Institut 
Fourier 53 (2003), 1387--1424.

\item{10.} Iftimie D.: The 3d Navier-Stokes equations seen as a perturbation 
of the 2d Navier-Stokes equations, Bull. Soc. Math. France 127 (1999), 
473--517.

\item{11.} Karch G. and Pilarczyk D.: Asymptotic stability of Landau 
solutions to Navier-Stokes system, Arch. Rational Mech. Anal. 202 (2011), 
115--131.

\item{12.} Krylov, N. V.: The Calderon-Zygmund theorem and its application 
to parabolic equations, Algebra i Analiz 13 (2001), 1--25 (in Russian).

\item{13.} Krylov, N. V.: The heat equation in $L_q(0,T;L_p(\Omega))$-spaces 
with weights, SIAM J. Math. Anal. 32 (2001), 1117--1141.

\item{14.} Ladyzhenskaya, O. A.; Solonnikov, V. A.; Uraltseva, N. N.: Linear 
and quasilinear equations of parabolic type, Nauka, Moscow 1967 (in Russian).

\item{15.} Mucha P. B.: Stability of nontrivial solutions of the Navier-Stokes 
system on the three-dimensional torus, J. Differential Equations 172 (2001), 
359--375.

\item{16.} Mucha P. B.: Stability of constant solutions to the Navier-Stokes 
system in $\R^3$, Appl. Math. 28 (2001), 301--310.

\item{17.} Mucha P. B.: Stability of 2d incompressible flows in $R^3$, J. 
Diff. Eqs. 245 (2008), 2355--2367.

\item{18.} Ponce G., Racke R., Sideris T. C. and Titi E. S.: Global stability 
of large solutions to the 3d Navier-Stokes equations, Comm. Math. Phys. 159 
(1994), 329--341.

\item{19.} Solonnikov, V. A.: Estimates of solutions of the Stokes equations 
in Sobolev spaces with a mixed norm, Zap. Nauchn. Sem. POMI 288 (2002), 
204--231.

\item{20.} Solonnikov, V. A.: On the estimates of solutions of nonstationary 
Stokes problem in anisotropic S.L. Sobolev spaces and on the estimate of 
resolvent of the Stokes problem, Usp. Mat. Nauk 58 (2) (2003) (350), 123--156.

\item{21.} Solonnikov, V. A.: Overdetermined elliptic boundary value 
problems, Zap. Nauchn. Sem. LOMI 21 (1971), 112--158 (in Russian).

\item{22.} Zadrzy\'nska, E.; Zaj\c aczkowski, W. M.: Stability of 
two-dimensional Navier-Stokes motions in the periodic case, J. Math. Anal. 
Appl. (2015), http://dx.doi.org/10.1016/j.jmaa.2014.10.026.

\item{23.} Zaj\c aczkowski, W. M.: Existence and regularity of solutions of 
some elliptic system in domains with edges, Diss. Math. 274 (1988), 91 pp.

\item{24.} Zaj\c aczkowski, W. M.: Some stability problem to the 
Navier-Stokes equations in the periodic case (to be published).

\item{25.} Zhou Y.: Asymptotic stability for the Navier-Stokes equations 
in the marginal class, Proc Roy. Soc. Edinburgh 136 (2006), 1099--1109.

\bye